\documentclass[leqno,12pt]{article}
\usepackage{}

\usepackage{amssymb}
\usepackage{euscript}
\usepackage[dvips]{graphicx}
\usepackage{flafter}
\usepackage{pstricks}
\usepackage[latin1]{inputenc}
\usepackage{amsmath}
\usepackage{amsfonts}
\usepackage{enumitem}
\usepackage{amssymb}
\usepackage{amsthm,array}

\usepackage{verbatim}

\usepackage{mathrsfs} 

\evensidemargin\oddsidemargin

\begin{document}

\baselineskip=18pt
\setcounter{page}{1}

\renewcommand{\theequation}{\thesection.\arabic{equation}}
\newtheorem{theorem}{Theorem}[section]
\newtheorem{lemma}[theorem]{Lemma}
\newtheorem{proposition}[theorem]{Proposition}
\newtheorem{corollary}[theorem]{Corollary}
\newtheorem{remark}[theorem]{Remark}
\newtheorem{fact}[theorem]{Fact}

\newcommand{\eqnsection}{
\renewcommand{\theequation}{\thesection.\arabic{equation}}
    \makeatletter
    \csname  @addtoreset\endcsname{equation}{section}
    \makeatother}
\eqnsection

\def\r{{\mathbb R}}
\def\e{{\mathbb E}}
\def\p{{\mathbb P}}
\def\q{{\mathbb Q}}
\def\f{{\mathcal F}}
\def\P{{\bf P}}
\def\E{{\bf E}}
\def\Q{{\bf Q}}
\def\z{{\mathbb Z}}
\def\N{{\mathbb N}}
\def\T{{\mathbb T}}
\def\root{{\varnothing}}
\def\S{\mathscr{S}}
\def\G{{\mathbb G}}
\def\L{{\mathbb L}}
\def\ev{{\mathscr E}}
\def\l{{\mathcal L}}

\def\deg{\chi}

\def\ee{\mathrm{e}}
\def\d{\, \mathrm{d}}



\vglue50pt

\centerline{\large\bf A necessary and sufficient condition for the non-trivial limit}
\centerline{\large\bf of the derivative martingale in a branching random walk}

\bigskip
\bigskip

\centerline{Xinxin Chen}

\medskip

\centerline{\it LPMA, Universit\'e Paris VI}

\bigskip
\bigskip
\bigskip

{\leftskip=2truecm \rightskip=2truecm \baselineskip=15pt \small

\noindent{\slshape\bfseries Summary.} We consider a branching random walk on the line. Biggins and Kyprianou \cite{biggins-kyprianou04} proved that, in the boundary case, the associated derivative martingale converges almost surly to a finite nonnegative limit, whose law serves as a fixed point of a smoothing transformation (Mandelbrot's cascade). In the present paper, we give a necessary and sufficient condition for the non-triviality of this limit and establish a Kesten-Stigum-like result.

\bigskip

\bigskip

}
\section{Introduction}
\label{s:intro}

$\phantom{aob}$ We consider a discrete-time branching random walk (BRW) on the real line, which can be described in the following way. An initial ancestor, called the root and denoted by $\root$, is created at the origin. It gives birth to some children which form the first generation and whose positions are given by a point process $\l$ on $\r$. For any integer $n\geq1$, each individual in the $n$th generation gives birth independently of all others to its own children in the $(n+1)$th generation, and the displacements of its children from this individual's position is given by an independent copy of $\l$. The system goes on if there is no extinction. We thus obtain a genealogical tree, denoted by $\T$. For each vertex (individual) $u\in\T$, we denote its generation by $|u|$ and its position by $V(u)$. In particular, $V(\root)=0$ and $(V(u); |u|=1)=\l$.

Note that the point process $\l$ plays the same role in the BRW as the offspring distribution in a Galton-Watson process. We introduce the Laplace-Stieltjes transform of $\l$ as follows:
\begin{equation}\label{def:lstransform}
\Phi(t):=\E\Big[\int_\r \ee^{-tx}\l(\d x)\Big]=\E\Big[\sum_{|u|=1}\ee^{-tV(u)}\Big], \text{ for } \forall t\in\r.
\end{equation}
Let $\Psi(t):=\log \Phi(t)$. We always assume in this paper $\Psi(0)>0$ so that $\E\Big[\sum_{|u|=1}1\Big]>1$. This yields that with strictly positive probability, the system survives. Let $q$ be the probability of extinction. Clearly, $q<1$.

Let $(\f_n; n\geq0)$ be the natural filtration of this branching random walk, i.e. let $\f_n:=\sigma\{(u, V(u)); |u|\leq n\}$. We introduce the additive martingale for any $t\in\r$,
\begin{equation}\label{def:additivemart}
W_n(t):=\sum_{|u|=n}\ee^{-t V(u)-n\Psi(t)}.
\end{equation}
It is a nonnegative martingale with respect to $(\f_n; n\geq0)$, which converges almost surely to a finite nonnegative limit. Biggins \cite{biggins77} established a necessary and sufficient condition for the mean convergence of $W_n(t)$, and generalized Kesten-Stigum theorem for the Galton-Watson processes. A simpler proof based on a change of measures was given later by Lyons \cite{lyons1997}.

More generally, Biggins and Kyprianou \cite{biggins-kyprianou04} studied the martingales produced by the so-called mean-harmonic functions. Given suitable conditions on the offspring distribution $\mathcal{L}$ of the branching random walk, like the $X\log X$ condition of the Kesten-Stigum theorem, they gave a general treatment to obtain the mean convergence of these martingales. In this paper, following their ideas, we work on one special example and give a Kesten-Stigum-like theorem.

Throughout this paper, we consider the boundary case (in the sense of \cite{biggins-kyprianou05}) where $\Psi(1)=\Psi^\prime(1)=0$, i.e.,
\begin{equation}\label{hyp:boundarycase}
\E\Big[\sum_{|u|=1}\ee^{-V(u)}\Big]=1, \qquad \E\Big[\sum_{|u|=1}V(u)\ee^{-V(u)}\Big]=0.
\end{equation}
In addition, we assume that
\begin{equation}\label{hyp:secondmom}
\sigma^2:=\E\Big[\sum_{|u|=1}V(u)^2\ee^{-V(u)}\Big]\in(0,\infty).
\end{equation}
We are interested in the derivative martingale, which is defined as follows:
\begin{equation}\label{def:derivativemart}
D_n:=\sum_{|u|=n} V(u)\ee^{-V(u)},\quad \forall n\geq0.
\end{equation}
It is a signed martingale with respect to $(\f_n)$, of mean zero. By Theorem 5.1 of \cite{biggins-kyprianou04}, under (\ref{hyp:boundarycase}) and (\ref{hyp:secondmom}), $D_n$ converges almost surely to a finite nonnegative limit, denoted by $D_\infty$. Moreover, $D_\infty$ satisfies the following equation (Mandelbrot's cascade):
\begin{equation}\label{eq:cascade}
D_\infty=\sum_{|u|=1} \ee^{-V(u)}D_\infty^{(u)},
\end{equation}
where $D_\infty^{(u)}$ are copies of $D_\infty$ independent of each other and of $\f_1$. Note that $D_\infty$ serves as a nonnegative fixed point of a smoothing transformation. From this point of view, the questions concerning the existence, uniqueness and asymptotic behavior of such fixed points have been much studied in the literature (\cite{biggins-kyprianou97,biggins-kyprianou05,liu98,liu00}). We are interested in the existence of a non-trivial fixed point, and we are going to determine when $\P(D_\infty>0)>0$.

It is known that $\P(D_\infty=0)$ is equal to either the extinction probability $q$ or 1 (see \cite{aidekon}, for example). We say that the limit $D_\infty$ is non-trivial if $\P(D_\infty>0)>0$, which means that $\P(D_\infty=0)=q$. Otherwise, it is trivially zero. In this paper, we give a sufficient and necessary condition for the non-triviality of $D_\infty$. The main result is stated as follows.

For any $y\in\r$, let $y_+:=\max\{y,0\}$ and let $\log_+y:=\log (\max\{y,1\})$. We introduce the following random variables:
\begin{equation}\label{def:rvs}
Y:=\sum_{|u|=1}\ee^{-V(u)},\quad Z:=\sum_{|u|=1}V(u)_+\ee^{-V(u)}.
\end{equation}

\begin{theorem}\label{thm:main}
The limit of the derivative martingale $D_n$ is non-trivial, namely $\P(D_\infty>0)>0$, if and only if the following condition holds:
\begin{equation}\label{cond:iff}
\E\Big(Z\log_+ Z+Y(\log_+ Y)^2\Big)<\infty.
\end{equation}
\end{theorem}

\begin{remark} In \cite{biggins-kyprianou04}, the authors studied the optimal condition for the non-triviality of $D_\infty$. However, there is a small gap between the necessary condition and the sufficient condition for $\P(D_\infty>0)>0$ in their Theorem 5.2. Our result fills this gap and gives the analogue of the result of \cite{Ren-Yang2011} in the case of branching Brownian motion.
\end{remark}

\begin{remark}
A\"id\'ekon proved that the condition (\ref{cond:iff}) is sufficient for $\P(D_\infty>0)>0$ (see Proposition A.3 in the Appendix of \cite{aidekon}).
\end{remark}

The paper is organized as follows. Section \ref{s:tool} introduces a change of measures based on a truncated martingale which is closely related to the derivative martingale. We also prove a proposition concerning certain behaviors of a centered random walk conditioned to stay positive at the end of Section \ref{s:tool}. Then, by using this proposition, we prove Theorem \ref{thm:main} in Section \ref{s:proofofthm}.

Throughout the paper, $(c_i)_{i\geq0}$ denote positive constants. We write $\E[f;\,A]$ for $\E[f1_{A}]$ and set $\sum_{\emptyset}:=0$.

\section{Lyons' change of measures via truncated martingales}
  \label{s:tool}

\subsection{Truncated martingales}

We begin with the well-known many-to-one lemma. For any $a\in\r$, let $\P_a$ be the probability measure such that $\P_a\Big(\big(V(u),\; u\in\T\big)\in\cdot\Big)=\P\Big(\big(V(u)+a,\; u\in\T\big)\in\cdot\Big)$. The corresponding expectation is denoted by $\E_a$. We write $\P$, $\E$ instead of $\P_0$, $\E_0$ for brevity. For any particle $u\in\T$, we denote by $u_i$ its ancestor at the $i$th generation, for $0\leq i<|u|$. In addition, we write $u_{|u|}:=u$. We thus denote its ancestral line by $[\![ \varnothing, \, u]\!]:=\{u_0,\; u_1,\cdots, u_{|u|}\}$.

\begin{lemma}[Many-to-one]\label{many-to-one}
There exists a sequence of i.i.d centered random variables $(S_{k+1}-S_{k})$, $k\geq 0$ such that for any $n\geq 1$ and any measurable function $g:\r^n\rightarrow \r_+$, we have
\begin{equation}
\E_a\Big[\sum_{|u|=n}g\big(V(u_1),\cdots, V(u_n)\big)\Big]=\E_a\Big[\ee^{S_n-a}g(S_1,\cdots, S_n)\Big],
\end{equation}
with $\P_a[S_0=a]=1$.
\end{lemma}

In view of (\ref{hyp:secondmom}), $S_1-S_0$ has a finite variance $\sigma^2=\E[S_1^2]=\E[\sum_{|u|=1}V(u)^2\ee^{-V(u)}]$.

Let $U^-(\d y)$ be  the renewal measure associated with the weak descending ladder height process of $(S_n,\; n\geq 0)$. Following the arguments in Section 2 of \cite{biggins03}, we obtain that for any measurable function $f:\r\rightarrow\r_+$,
\begin{equation}\label{eq:renew}
\E\Big[\sum_{j=0}^{\tau-1}f(-S_j)\Big]=\int_0^\infty f(y)U^-(\d y),
\end{equation}
where $\tau$ be the first time that $(S_n)$ enters $(0,\infty)$, namely $\tau:=\inf\{k> 0,\; S_k\in(0,\infty)\}$ which is proper here.
We define $R(x):=U^-([0,x))$ for all $x>0$ and define $R(0):=1$. Note that $R(x)$ equals the renewal function $U^-([0,x])$ at points of continuity. We collect the following properties of this function $R(x)$ which are consequences of the renewal theorem (see \cite{biggins03, bertoin-doney, tanaka}).
\begin{fact}\label{fact:base}
\begin{itemize}[fullwidth]
  \item [(i)] There exists a positive constant $c_0>0$ such that
  \begin{equation}
  \lim_{x\rightarrow\infty}\frac{R(x)}{x}=c_0.
  \end{equation}
  \item [(ii)] There exist two constants $0<c_1<c_2<\infty$ such that
  \begin{equation}\label{eq:upforr}
  c_1(1+x)\leq R(x)\leq c_2(1+x),\qquad \forall x\geq 0.
  \end{equation}
  \item [(iii)] For any $x\geq 0$, we have $\E[R(S_1+x)1_{(S_1+x>0)}]=R(x)$.
\end{itemize}
\end{fact}

Let $\beta\geq 0$. Started from $V(\root)=a$, we add a barrier at $-\beta$ to the branching random walk. Now, we define the following truncated random variables:
\begin{equation}
D_n^{(\beta)}:=\sum_{|x|=n}R(V(x)+\beta)\ee^{-V(x)}1_{(\min_{1\leq k\leq n}V(x_k)>-\beta)},\qquad \forall n\geq 1,
\end{equation}
and $D_0^{(\beta)}:=R(a+\beta)\ee^{-a}1_{(a\geq-\beta)}$.

\begin{lemma}
For any $a\geq0$ and $\beta\geq0$, under $\P_a$, the process $(D_n^{(\beta)},\; n\geq 0)$ is a nonnegative martingale with respect to $(\f_n,\; n\geq 0)$.
\end{lemma}

This lemma follows immediately from (iii) of Fact \ref{fact:base} and the branching property. We feel free to omit its proof and call $(D_n^{(\beta)})$ the truncated martingale. It also tells us that under $\P_a$, $(D_n^{(\beta)},\; n\geq 0)$ converges almost surely to a finite nonnegative limit, which we denote by $D_\infty^{(\beta)}$.

The connection between the limits of the derivative martingale and truncated martingales is recorded in the following Lemma, the proof of which can be referred to \cite{biggins-kyprianou04} and \cite{aidekon}.

\begin{lemma}
\begin{itemize}
  \item [(1)] If $D_\infty$ is trivial, i.e., $\P(D_\infty=0)=1$, then for any $\beta\geq 0$, $D_\infty^{(\beta)}$ is trivially zero under $\P$.
  \item [(2)] Under $\P$, if there exists some $\beta\geq 0$ such that $D_\infty^{(\beta)}$ is trivially zero, so is $D_\infty$.
\end{itemize}
\label{lem:connection}
\end{lemma}

Thanks to Lemma \ref{lem:connection}, we only need to investigate the truncated martingale $(D_n^{(0)}; n\geq0)$ and determine when its limit is non-trivial.

\subsection{Lyons' change of probabilities and spinal decomposition}

Let $\beta=0$. With this nonnegative martingale $(D_n^{(0)},\; n\geq 0)$, we define for any $a\geq0$ a new probability measure $\Q_a$ such that for any $n\geq 1$,
\begin{equation}
\frac{d\Q_a}{d\P_a}\Big\vert_{\f_n}=\frac{D_n^{(0)}}{R(a)e^{-a}}.
\end{equation}

$\Q_a$ is defined on $\f_\infty:=\vee_{n\geq 0}\f_n$. Let us give an intuitive description of the branching random walk under $\Q_a$, which is known as the spinal decomposition. We start from one single particle $\omega_0$, located at the position $V(\omega_0)=a$. At time 1, it dies and produces a point process distributed as $(V(u);\;|u|=1)$ under $\Q_a$. Among the children of $\omega_0$, $\omega_1$ is chosen to be $u$ with probability proportional to $R(V(u))\ee^{-V(u)}1_{(V(u)>0)}$. At each time $n+1$, each particle $v$ in the $n$th generation dies and produces independently a point process distributed as $(V(u);\; |u|=1)$ under $\P_{V(v)}$ except $\omega_n$, which dies and generates independently a point process distributed as $(V(u);\;|u|=1)$ under $\Q_{V(\omega_n)}$. And then $\omega_{n+1}$ is chosen to be $u$ among the children of $\omega_n$, with probability proportional to $R(V(u))\ee^{-V(u)}1_{(\min_{1\leq k\leq n+1}V(u_k)>0)}$. We still use $\T$ to denote the genealogical tree. Then $(\omega_n;\; n\geq0)$ is an infinite ray in $\T$, which is called the spine. The rigorous proof was given in Appendix A of \cite{aidekon}. Indeed, this type of measures' change and the establishment of a spinal decomposition have been developed in various cases of the branching framework; see, for example \cite{lyons1997, hu-shi, chauvin-rouault, roberts-harris}.

We state the following fact about the distribution of the spine process $(V(\omega_n); \; n\geq 0)$ under $\Q_a$.

\begin{fact}
Let $a\geq 0$. For any $n\geq 0$ and any measurable function $g:\r^{n+1}\rightarrow\r_+$, we have
\begin{equation}\label{rwstaypos}
\E_{\Q_a}\Big[g(V(\omega_0),\cdots,V(\omega_n))\Big]=\frac{1}{R(a)}\E_a\Big[g(S_0,\cdots, S_n)R(S_n);\; \min_{1\leq k\leq n}S_i>0\Big],
\end{equation}
where $(S_n)$ is the same as that in Lemma \ref{many-to-one}.
\end{fact}
For convenience, let $(\zeta_n;\; n\geq 0)$ be a stochastic process under $\P_a$ such that
\begin{equation}\label{eq:simple}
\P_a[\Big(\zeta_n;\; n\geq 0\Big)\in \cdot]=\Q_a[\Big(V(\omega_n);\; n\geq 0 \Big)\in \cdot].
\end{equation}
Obviously, under $\P_a$, $(\zeta_n;\;n\geq0)$ is a Markov chain with transition probabilities $P$ so that, for any $x\geq 0$, $P(x,\d y)=\frac{R(y)}{R(x)}1_{(y>0)}\P_x(S_1\in \d y)$. This process $(\zeta_n)$ is usually called a random walk conditioned to stay positive. It has been arisen and studied in, for instance, \cite{tanaka, bertoin-doney, biggins03, vatutin-wachtel}. In what follows, we state some results about $(\zeta_n)$, which will be useful later in Section \ref{s:proofofthm}.

\subsection{Random walk conditioned to stay positive}

Recall that $(S_n)$ is a centered random walk on $\r$ with finite variance $\sigma^2$. Let $\tau_-$ be the first time that $(S_n)$ hits $(-\infty,0]$, namely, $\tau_-:=\inf\{k\geq1:\; S_k\leq0\}$. Let $(T_k,\, H_k;\; k\geq0)$ be the strict ascending ladder epochs and heights of $(S_n;\;n\geq0)$, i.e., $T_0=0$, $H_0:=S_0$ and for any $k\geq1$, $T_k:=\inf\{j>T_{k-1}:\, S_j>H_{k-1}\}$, $H_k:=S_{T_k}$. We denote by $U(\d x)$ the corresponding renewal measure (see Chapter XII in \cite{feller1971}, for example). Then, similarly to (\ref{eq:renew}), for any measurable function $f:\r\rightarrow\r_+$,
\begin{equation}\label{eq:renewal}
\E\Big[\sum_{n=0}^{\tau_--1}f(S_n)\Big]=\E\Big[\sum_{k\geq0}f(H_k)\Big]=\int_0^\infty f(x)U(\d x).
\end{equation}
We deduce from (\ref{rwstaypos}) and (\ref{eq:renewal}) that
\begin{align}\label{eq:expsum}
\E\Big[\sum_{n\geq 0}f(\zeta_n)\Big]&=\E_{\Q_0}\Big[\sum_{n\geq 0}f(V(\omega_n))\Big]=\sum_{n\geq 0}\E\Big[f(S_n)R(S_n)1_{(\min_{1\leq k\leq n}S_k>0)}\Big]\nonumber\\
                                    &= \E\Big[\sum_{n=0}^{\tau_--1}f(S_n)R(S_n)\Big]= \int_0^\infty f(x)R(x)U(\d x).
\end{align}

Recall also that $U^-(\d x)$ is the renewal measure associated with the weak descending ladder height process of $(S_n)$. By the renewal theorem (see P.360 in \cite{feller1971}), there exist two constants $c_3$, $c_4>0$ such that for $\forall x,\, y\geq0$,
\begin{eqnarray}
c_3(1+x)\leq U([0,x])\leq c_4(1+x), &\quad& 0\leq U([x,x+y])\leq c_4(1+y);\label{eq:upforu}\\
c_3(1+x)\leq U^-([0,x])\leq c_4(1+x), &\quad& 0\leq U^-([x,x+y])\leq c_4(1+y).\label{eq:upforu-}
\end{eqnarray}

Given a non-increasing function $F\geq0$, we present the following proposition, which gives a necessary and sufficient condition for the infinity of the series $\sum_n F(\zeta_n)$.

\begin{proposition}\label{prop:spine}
Let $F: [0,\infty)\rightarrow[0,\infty)$ be non-increasing. Then
\begin{equation}\label{eq:eqiv}
\int_0^\infty F(y)y\d y=\infty \Longleftrightarrow \sum_{n\geq 0}F(\zeta_n)=\infty,\quad \P\textrm{-a.s.}
\end{equation}
\end{proposition}

Note that $(\zeta_n)$ can be viewed as a discrete-time counterpart of the 3-dimensional Bessel process, for which a similar result holds (see, for instance, Ex 2.5, Chapter XI of \cite{revuz-yor}). And we will prove (\ref{eq:eqiv}) in a similar way as for the Bessel process.

\begin{proof} Observe that $0\leq F(x)\leq F(0)<\infty$ for any $x\geq0$. So there is no difference between the two events $\{\sum_{n\geq0}F(\zeta_n)=\infty\}$ and $\{\sum_{n\geq1}F(\zeta_n)=\infty\}$.

We first prove ``$\Longleftarrow$" in (\ref{eq:eqiv}).
It follows from (\ref{eq:upforr}) and (\ref{eq:upforu}) that
\begin{equation}
\int_0^\infty F(y)y\d y=\infty \Longleftrightarrow \int_0^\infty F(y)R(y)U(\d y)=\infty.
\end{equation}
Actually, by (\ref{eq:expsum}),
\begin{equation*}
\E\Big[\sum_{n\geq0}F(\zeta_n)\Big]=\int_0^\infty F(y)R(y)U(\d y).
\end{equation*}
Clearly, $\P\Big[\sum_{n\geq 0}F(\zeta_n)=\infty\Big]=1$ yields $\int_0^\infty F(y)R(y)U(\d y)=\infty$. The ``$\Longleftarrow$" in (\ref{eq:eqiv}) is hence proved.

To prove ``$\Longrightarrow$" in (\ref{eq:eqiv}), we only need to show that if $\P\Big[\sum_{n\geq 0}F(\zeta_n)=\infty\Big]<1$, then $\int_0^\infty F(y)y\d y<\infty$. From now on, we suppose that $\P\Big[\sum_{n\geq 0}F(\zeta_n)=\infty\Big]<1$, which is equivalent to say that,
\begin{equation}\label{eq:iffini}
\P\bigg[\sum_{n\geq 1}F(\zeta_n)<\infty\bigg]>0.
\end{equation}

We draw support from Tanaka's construction for the random walk conditioned to stay positive (\cite{tanaka, biggins03}). Recall that $\tau=\inf\{k\geq1: S_k\in(0,\infty)\}$. We hence obtain an excursion $(S_j;\; 0\leq j\leq \tau)$, which is denoted by $\xi=(\xi(j),\; 0\leq j\leq \tau)$. Let $\{\xi_k=(\xi_k(j),\,0\leq j\leq\tau_k);\;k\geq1\}$ be a sequence of independent copies of $\xi$. For any $k\geq1$, let
\begin{equation}\label{eq:reflet}
\nu_k(j):=\xi_k(\tau_k)-\xi_k(\tau_k-j), \ \forall 0\leq j\leq \tau_k.
\end{equation}
This brings out another sequence of i.i.d. excursions $\{\nu_k=(\nu_k(j),\; 0\leq j\leq \tau_k);\, k\geq1\}$, based on which we reconstruct the random walk conditioned to stay position $(\zeta_n)$ in the following way. Define for any $k\geq1$,
\begin{align}
T_k^+&:=\tau_1+...+\tau_k;\\
H_k^+&:=\nu_1(\tau_1)+...+\nu_k(\tau_k)=\xi_1(\tau_1)+...+\xi_k(\tau_k),
\end{align}
and let $T_0^+=H_0^+=0$. Then the process
\begin{equation}
\zeta_n=H_k^++\nu_{k+1}(n-T_k^+), \quad \textrm{for }T_k^+<n\leq T_{k+1}^+,
\end{equation}
with $\zeta_0=0$, is what we need.

We actually establish un process distributed as $(\zeta_n)$. For brevity, we still denote it by $(\zeta_n)$ without changing any conclusion in this proof. For any $k\geq1$, let
\begin{equation}
\chi_k(F):=\sum_{n=T_{k-1}^++1}^{T_{k}^+}F(\zeta_n)=\sum_{j=1}^{\tau_{k}}F\Big(H_{k-1}^++\nu_{k}(j)\Big),
\end{equation}
so that $\sum_{n\geq1}F(\zeta_n)=\sum_{k\geq1}\chi_k(F)$.

By (\ref{eq:reflet}), we get that
\begin{eqnarray*}
\chi_k(F)&=&\sum_{j=1}^{\tau_{k}}F\Big(H_{k-1}^++\xi_{k}(\tau_{k})-\xi_{k}(\tau_{k}-j)\Big)\\
&=&\sum_{j=0}^{\tau_{k}-1}F\Big(H_{k}^+-\xi_{k}(j)\Big).
\end{eqnarray*}

(\ref{eq:iffini}) hence becomes that
\begin{equation}
\P\bigg[\sum_{k\geq1}\chi_k(F)<\infty\bigg]=\P\bigg[\sum_{k\geq1}\sum_{j=0}^{\tau_k-1}F\Big(H_k^+-\xi_k(j)\Big)<\infty\bigg]>0.
\end{equation}
By Theorem 1 in Chapter XVIII.5 of \cite{feller1971}, as $(S_n)$ is of finite variance, we have $b^+:=\E[H_1^+]=\E[S_\tau]<\infty$. It follows from Strong Law of Large Numbers that $\P$-a.s.,
\begin{equation}
\lim_{k\rightarrow\infty}\frac{H_k^+}{k}= b^+.
\end{equation}
Let $A>\max\{1,b^+\}$. This tells us that $\P$-a.s., for all large $k$, $H_k^+\leq Ak$. As $F$ is non-increasing, one sees that
\begin{equation}
\P\bigg[\sum_{k\geq1}\sum_{j=0}^{\tau_k-1}F\Big(Ak-\xi_k(j)\Big)<\infty\bigg]\geq \P\bigg[\sum_{k\geq1}\sum_{j=0}^{\tau_k-1}F\Big(H_k^+-\xi_k(j)\Big)<\infty\bigg]>0.
\end{equation}
For any $k\geq1$ let
\begin{equation}
\widetilde{\chi}_k:=\sum_{j=0}^{\tau_k-1}F\Big(Ak-\xi_k(j)\Big).
\end{equation}
So, $\P\Big[\sum_{k\geq1}\widetilde{\chi}_k<\infty\Big]>0$. Recall that $\{\xi_k, k\geq1\}$ is a sequence of independent copies of $(S_j;\; 0\leq j\leq\tau)$. This yields the independence of the sequence $\{\widetilde{\chi}_k, k\geq1\}$. It follows from Kolmogorov's $0$-$1$ law that
\begin{equation}\label{newassum}
\P\bigg[\sum_{k\geq1}\sum_{j=0}^{\tau_k-1}F(Ak-\xi_k(j))<\infty\bigg]=\P\bigg[\sum_{k\geq1}\widetilde{\chi}_k<\infty\bigg]=1.
\end{equation}

Moreover, let $E_M:=\Big\{\sum_{k\geq1}\widetilde{\chi}_k< M\Big\}$ for any $M>0$. Either there exists some $M_0<\infty$ such that $\P[E_{M_0}]=1$, or $\P[E_M]<1$ for all $M\in(0,\infty)$. On the one hand, if $\P[E_{M_0}]=1$ for some $M_0<\infty$, then
\begin{eqnarray*}
M_0&\geq&\E\bigg[\sum_{k\geq1}\widetilde{\chi}_k\bigg]=\E\bigg[\sum_{k\geq1}\sum_{j=0}^{\tau_k-1}F\Big(Ak-\xi_k(j)\Big)\bigg]\\
&=&\sum_{k\geq1}\E\bigg[\sum_{j=0}^{\tau-1}F(Ak-S_j)\bigg]\\
&=&\sum_{k\geq1}\int_0^\infty F(Ak+y)U^-(dy),
\end{eqnarray*}
where the last equality follows from (\ref{eq:renewal}). One sees that $\sum_{k\geq1}\int_0^\infty F(Ak+y)U^-(\d y)<\infty$. It follows from the renewal theorem that there exists $B>0$ such that $U^-([jB,jB+B))>\delta>0$ for any $j\geq0$. As $F$ is non-increasing,
 \begin{equation}
 \sum_{k\geq1}\sum_{j\geq1}F(Ak+Bj)\delta\leq \sum_{k\geq1}\int_0^\infty F(Ak+y)U^-(\d y)<\infty.
 \end{equation}
We hence observe that $\int_{A}^\infty \d z\int_{B}^\infty F(y+z)\d y\leq \sum_{k\geq1}\sum_{j\geq1}F(Ak+Bj)AB<\infty$. This implies that
$$\int_0^\infty F(x)x\d x=\int_0^\infty dz \int_0^\infty F(z+y)\d y\leq F(0)AB+\int_A^\infty \d z\int_B^\infty F(y+z)\d y<\infty,$$
which is what we need.

On the other hand, if $\P[E_M]<1$ for all $M\in(0,\infty)$, we have $\lim_{M\uparrow\infty}\P[E_M]=1$ because of (\ref{newassum}). For any $k\geq1$ and any $\ell\geq1$, define:
\begin{equation}
\Lambda_\ell^{(k)}:=\sum_{j=0}^{\tau_k-1}1_{\{A(\ell-1)\leq-\xi_k(j)<A\ell\}}.
\end{equation}
As $\sum_{\ell\geq1} 1_{\{A(\ell-1)\leq-\xi_k(j)<A\ell\}}=1$, we get that for any $k\geq1$,
\begin{eqnarray*}
\widetilde{\chi}_k&=&\sum_{j=0}^{\tau_k-1}F\Big(Ak-\xi_k(j)\Big)\sum_{\ell\geq1} 1_{\{A(\ell-1)\leq-\xi_k(j)<A\ell\}}\\
&=&\sum_{\ell\geq1} \sum_{j=0}^{\tau_k-1}F\Big(Ak-\xi_k(j)\Big) 1_{\{A(\ell-1)\leq-\xi_k(j)<A\ell\}}\\
&\geq&\sum_{\ell\geq1} F(Ak+A\ell) \Lambda_\ell^{(k)},
\end{eqnarray*}
where the last inequality holds because $F$ is non-increasing. It follows that
\begin{eqnarray}
\sum_{k\geq1}\widetilde{\chi}_k &\geq& \sum_{k\geq1}\sum_{\ell\geq1} F(Ak+A\ell) \Lambda_\ell^{(k)}= \sum_{n=2}^\infty F(An)\sum_{k=1}^{n-1}\Lambda_{n-k}^{(k)}\nonumber\\
&=& \sum_{m=1}^\infty F(Am+A)m Y_m,\label{eq:lowerbound}
\end{eqnarray}
where
\begin{equation}
Y_m:=\frac{\sum_{k=1}^m \Lambda_{m+1-k}^{(k)}}{m},\  \forall m\geq1.
\end{equation}
We claim that there exists a $M_1>0$ sufficiently large such that for any $m\geq1$,
\begin{equation}\label{eq:uniformcontrol}
c_6\geq\E\Big[Y_m1_{E_{M_1}}\Big]\geq c_5>0,
\end{equation}
where $c_5$, $c_6$ are positive constants. We postpone the proof of (\ref{eq:uniformcontrol}) and go back to (\ref{eq:lowerbound}). It follows that
\begin{eqnarray}
M_1&\geq&\E\bigg[1_{E_{M_1}}\sum_{k\geq1}\widetilde{\chi}_k\bigg]\geq \E\bigg[1_{E_{M_1}}\sum_{m=1}^\infty F(Am+A)m Y_m\bigg]\nonumber\\
&\geq& \sum_{m\geq1}F(Am+A)m\E\Big[Y_m1_{E_{M_1}}\Big].
\end{eqnarray}
By (\ref{eq:uniformcontrol}), we obtain that
\begin{equation}
 \sum_{m\geq1}F(Am+A)m \leq M_1/c_5<\infty .
\end{equation}
This implies that $\int_0^\infty F(y)y\d y<\infty$ thus completes the proof of Proposition \ref{prop:spine}.

It remains to prove (\ref{eq:uniformcontrol}).

We begin with the first and second moments of $Y_m$. Since $\{\omega_k; k\geq1\}$ are i.i.d. copies of $(S_j, 0\leq j\leq \tau)$, $(\Lambda_\ell^{(k)}; \ell\geq1)$, $k\geq1$ are i.i.d. This yields that
\begin{eqnarray}
\E[Y_m]&=&\frac{1}{m}\sum_{k=1}^m \E[\Lambda_{m+1-k}^{(k)}]=\frac{1}{m}\sum_{k=1}^m \E[\Lambda_{m+1-k}^{(1)}]\nonumber\\
&=& \frac{1}{m}\E\bigg[\sum_{k=1}^m \Lambda_{k}^{(1)}\bigg]=\frac{1}{m}\E\bigg[\sum_{j=0}^{\tau-1}1_{\{-S_j<Am\}}\bigg]\nonumber\\
&=& \frac{R(Am)}{m}.
\end{eqnarray}
where the last equality comes from (\ref{eq:renew}). By (\ref{eq:upforr}), for any $m\geq1$,
\begin{equation}\label{eq:meanofY}
c_1 A \leq \E[Y_m]\leq c_2(A+1)=:c_6.
\end{equation}
Obviously, we have $\E[Y_m1_{E_M}]\leq c_6$ for any $m\geq1$ and any $M>0$. The fact that $\Lambda_\cdot^{(k)}, k\geq1$, are i.i.d. yields also that
\begin{equation}\label{varofY}
\textrm{Var}(Y_m)=\frac{1}{m^2}\sum_{k=1}^m \textrm{Var}(\Lambda_k^{(1)})\leq \frac{1}{m^2}\sum_{k=1}^m \E[(\Lambda_k^{(1)})^2].
\end{equation}

Note that $\Lambda_1^{(1)}$ is distributed as $\sum_{j=0}^{\tau-1}1_{\{-S_j<A\}}$ with $\tau=\inf\{k>0: S_k>0\}$. We see that
\begin{eqnarray}
\E\Big[\Big(\Lambda_1^{(1)}\Big)^2\Big]&=& \E\Big[\Big(\sum_{j=0}^{\tau-1}1_{\{-S_j<A\}}\Big)^2\Big]\nonumber\\
&\leq& 2\E\Big[\sum_{j=0}^{\tau-1}1_{\{-S_j<A\}}\sum_{k=j}^{\tau-1}1_{\{-S_k<A\}}\Big].\nonumber
\end{eqnarray}
By Markov property, we obtain that
\begin{equation}
\E\Big[\Big(\Lambda_1^{(1)}\Big)^2\Big]\leq 2\E\Big[\sum_{j=0}^{\tau-1}1_{\{-S_j<A\}} R(A, -S_j)\Big],\label{secmom}
\end{equation}
where
\begin{equation}
R(x, y):=\E\Big[\sum_{i=0}^{\tau_y-1}1_{\{S_i>y-x\}}\Big]\textrm{ with } \tau_y:=\inf\{k>0:\, S_k>y\}\textrm{ for } x,y\geq0.
\end{equation}
It follows from (\ref{eq:renew}) that
\begin{equation}\label{eq:secondmom}
\E\Big[\Big(\Lambda_1^{(1)}\Big)^2\Big]\leq 2 \int_0^A R(A, y)U^-(\d y).
\end{equation}
Consider now the strict ascending ladder epochs and heights $(T_k,\, H_k)$ of $(S_n)$. We get that
\begin{equation*}
R(x,y)=\E\Big[\sum_{k=0}^\infty 1_{\{y\geq H_k>y-x\}}\sum_{n=T_k}^{T_{k+1}-1}1_{\{S_n>y-x\}}\Big].
\end{equation*}
By applying the Markov property at the times $(T_k;\; k\geq1)$ and (\ref{eq:renew}), we have for $x,\, y\geq0$,
\begin{eqnarray}\label{eq:doublerenewal}
R(x,y)=\E\Big[\sum_{k\geq0}R(H_k+x-y)1_{\{y\geq H_k> y-x\}}\Big]=\int_{(y-x)_+}^y R(x-y+z) U(\d z).
\end{eqnarray}
Plugging it into (\ref{eq:secondmom}) then using (\ref{eq:upforr}), (\ref{eq:upforu-}) and (\ref{eq:upforu}) implies that
\begin{equation}
\E\Big[\Big(\Lambda_1^{(1)}\Big)^2\Big]\leq c_7(1+A)^3\leq c_8 A^3,
\end{equation}
(see also Lemma 2 in \cite{biggins03}).

Moreover, for any $\ell\geq2$, $\Lambda_{\ell}^{(1)}$ has the same law as $\sum_{j=0}^{\tau-1}1_{\{\ell A-A\leq -S_j<\ell A\}}$. Similarly, we get that
\begin{eqnarray*}
\E\Big[\Big(\Lambda_{\ell}^{(1)}\Big)^2\Big]&=&\E\Big[\Big(\sum_{j=0}^{\tau-1}1_{\{\ell A-A\leq-S_j<\ell A\}}\Big)^2\Big]\\
&\leq & 2\E\Big[\sum_{j=0}^{\tau-1}1_{\{\ell A-A\leq -S_j<\ell A\}}\sum_{k=j}^{\tau-1}1_{\{\ell A-A\leq -S_k<\ell A\}}\Big].
\end{eqnarray*}
Once again, by Markov property then by (\ref{eq:renew}),
\begin{eqnarray*}
\E\Big[\Big(\Lambda_{\ell}^{(1)}\Big)^2\Big]&\leq&2\E\Big[\sum_{j=0}^{\tau-1}1_{\{\ell A-A\leq -S_j<\ell A\}}\Big(R(\ell A,-S_j)-R(\ell A-A, -S_j)\Big)\Big]\\
&=& 2\int_{\ell A- A}^{\ell A}\Big(R(\ell A,y)-R(\ell A-A, y)\Big)U^-(\d y).
\end{eqnarray*}
Plugging (\ref{eq:doublerenewal}) into it yields that for $\ell\geq2$,
\begin{eqnarray*}
\E\Big[\Big(\Lambda_{\ell}^{(1)}\Big)^2\Big]&\leq&2\int_{\ell A-A}^{\ell A}\Big(\int_0^y R(\ell A-y+z)U(\d z)-\int_{y-\ell A+A}^y R(\ell A-A-y+z)U(\d z)\Big)U^-(\d y)\\
&=& 2\int_{\ell A-A}^{\ell A}\Big(\int_0^{y-\ell A+A} R(\ell A-y+z)U(\d z)\\
& &\quad \quad +\int_{y-\ell A+A}^y U^-([\ell A-A-y+z,\ell A-y+z))U(\d z)\Big)U^-(\d y),
\end{eqnarray*}
where the last equality holds because $R(x)=U^-([0,x))$. Observe that $R(\ell A-y+z)\leq R(A)$ for $0\leq z\leq y-\ell A+A$ and $\ell A-A\leq y\leq \ell A$. Recall that $A\geq1$. By (\ref{eq:upforr}), (\ref{eq:upforu}) and (\ref{eq:upforu-}),
\begin{eqnarray*}
\E\Big[\Big(\Lambda_{\ell}^{(1)}\Big)^2\Big]&\leq& c_9\int_{\ell A-A}^{\ell A} \Big(\int_0^{y-\ell A+A}(A+1)U(\d z)+\int_{y-\ell A+A}^{y}(1+ A) U(\d z)\Big) U^-(\d y)\\
&\leq &c_{10}(A+1)\int_{\ell A-A}^{\ell A} \Big(y+1\Big) U^-(\d y)\\
&\leq &c_{11} \ell A^3,
\end{eqnarray*}
with $c_{11}\geq c_8$. Going back to (\ref{varofY}), for any $m\geq1$,
\begin{equation}
\textrm{Var}(Y_m)\leq \frac{\sum_{\ell=1}^m c_{11}\ell A^3}{m^2}\leq c_{12} A^3.
\end{equation}
Combining this with (\ref{eq:meanofY}) implies that $\E[Y_m^2]= \textrm{Var}(Y_m)+\E[Y_m]^2\leq c_2^2(1+A)^2+ c_{12} A^3$. We then use Paley-Zygmund inequality to obtain that
\begin{equation}
\P\Big[Y_m > \frac{1}{2}\E[Y_m]\Big]\geq\frac{\E[Y_m]^2}{4\E[Y_m^2]}\geq\frac{c_1^2 A^2}{4(c_2^2(1+A)^2+ c_{12} A^3)}:=c_{13}>0.
\end{equation}
So for any $0\leq u \leq c_1 A/2\leq \E[Y_m]/2$, we have
\begin{equation}
\P\Big(Y_m\leq u \Big)\leq \P\Big(Y_m\leq \E[Y_m]/2 \Big)\leq 1-c_{13}.
\end{equation}
There exists $M_1>0$ such that $\P(E_{M_1})\geq 1-c_{13}/2$, since $\lim_{M\uparrow\infty}\P[E_M]=1$. For such $M_1>0$,
\begin{equation}
\E[Y_m1_{E_{M_1}}]=\E\Big[\int_0^{Y_m}1_{E_{M_1}}du\Big]=\int_0^\infty\P[\{Y_m >u\}\cap E_{M_1}]du.
\end{equation}
Notice that $\P[\{Y_m>u\}\cap E_{M_1}]\geq\Big(\P[E_{M_1}]-\P[Y_m\leq u]\Big)_+$, which is larger than $c_{13}/2$ when $0\leq u\leq c_1A/2$. As a consequence,
\begin{equation}
\E[Y_m1_{E_{M_1}}]=\int_0^\infty\P[\{Y_m >u\}\cap E_{M_1}]du \geq \int_0^{c_1 A/2} \frac{c_{13}}{2}du=\frac{c_1 c_{13} A}{4}=:c_5>0.
\end{equation}
This completes the proof of (\ref{eq:uniformcontrol}), hence completes the proof of ``$\Longrightarrow$" in (\ref{eq:eqiv}). Proposition \ref{prop:spine} is proved.
\end{proof}

\section{Proof of the main theorem}
\label{s:proofofthm}
Recall that we are in the regime that
\begin{equation}\label{eq:settings}
\E\Big[\sum_{|u|=1}\ee^{-V(u)}\Big]=1,
\quad
\E\Big[\sum_{|u|=1}V(u)\ee^{-V(u)}\Big]=0,
\quad
\sigma^2=\E\Big[\sum_{|u|=1}V(u)^2\ee^{-V(u)}\Big]<\infty.
\end{equation}
Recall also that equivalence in Theorem \ref{thm:main} is as follows:
\begin{equation}\label{eq:iff}
\E\Big[Y\Big(\log_+Y\Big)^2\Big]+\E\Big[Z\log_+Z\Big]<\infty\Longleftrightarrow \P[D_\infty>0]>0.
\end{equation}
with $Y=\sum_{|u|=1}\ee^{-V(u)}$ and $ Z=\sum_{|u|=1}V(u)_+\ee^{-V(u)}$.

This section is devoted to proving that the condition on the left-hand side of (\ref{eq:iff}) (i.e. (\ref{cond:iff})) is necessary and sufficient for mean convergence of the truncated martingale $\Big\{D_n^{(0)}=\sum_{|u|=n}R(V(u))\ee^{-V(u)}1_{\{V(u_k)>0,\forall 1\leq k\leq n\}};n\geq0\Big\}$. In view of Lemma \ref{lem:connection}, this follows the non-triviality of $D_\infty$, hence proves Theorem \ref{thm:main}.

In what follows, we state a result about the mean convergence of the truncated martingale $\Big\{D_n^{(0)}; n\geq0\Big\}$, which is one special case of Theorem 2.1 in Biggins and Kyprianou \cite{biggins-kyprianou04}.

Define
\begin{equation}
X:=\frac{D_1^{(0)}}{D_0^{(0)}}1_{(D_0^{(0)}>0)}+1_{(D_0^{(0)}=0)}.
\end{equation}
Then for any $a\geq0$, under $\P_a$, $X=\frac{\sum_{|u|=1}R(V(u))\ee^{-V(u)}1_{(V(u)>0)}}{R(a)\ee^{-a}}$. 

\begin{theorem}[Biggins and Kyprianou \cite{biggins-kyprianou04}]
\label{thm:BK}
$(\zeta_n)$ is a random walk conditioned to stay positive, whose law was given in (\ref{eq:simple}).
\begin{itemize}
\item[(i)] If \begin{equation}\label{eq:meanconv}
\P\textrm{-a.s. } \sum_{n\geq1}\E_{\zeta_n}\Big[X\Big(R(\zeta_n)\ee^{-\zeta_n}X\wedge 1\Big)\Big]<\infty,
\end{equation}
then $\E[D_\infty^{(0)}]=R(0)$.
\item[(ii)] If for all $y>0$,
\begin{equation}\label{intermedia}
 \P\textrm{-a.s. }\sum_{n=1}^\infty\E_{\zeta_n}\Big[X;\; R(\zeta_n)\ee^{-\zeta_n}X\geq y\Big]=\infty,
\end{equation}
then $\E[D_\infty^{(0)}]=0$.
\end{itemize}
\end{theorem}

Our proof relies on this theorem. First, in Subsection \ref{ss:proofofif}, we give a short proof for the sufficient part to accomplish our arguments even though it has already been proved in \cite{aidekon}. In Subsection \ref{ss:proofofonlyif}, we prove that (\ref{cond:iff}) is also the necessary condition by using Proposition \ref{prop:spine}.

\subsection{(\ref{cond:iff}) is a sufficient condition}
\label{ss:proofofif}
This subsection is devoted to proving that
\begin{equation}\label{eq:if}
\E\Big[Y\Big(\log_+Y\Big)^2\Big]+\E\Big[Z\log_+Z\Big]<\infty\Longrightarrow \E[D_\infty^{(0)}]=R(0)=1.
\end{equation}

\begin{proof}[Proof of (\ref{eq:if})]
According to (i) of Theorem \ref{thm:BK}, it suffices to show that
\begin{equation}\label{eq:finimean}
\E\Big[Y\Big(\log_+Y\Big)^2\Big]+\E\Big[Z\log_+Z\Big]<\infty\Longrightarrow \P\textrm{-a.s. } \sum_{n\geq0}\E_{\zeta_n}\Big[X\Big(R(\zeta_n)\ee^{-\zeta_n}X\wedge 1\Big)\Big]<\infty.
\end{equation}

 For any particle $x\in\mathbb{T}\setminus\{\varnothing\}$, we denote its parent by $\overleftarrow{u}$ and define its relative displacement by
\begin{equation}
\Delta V(u):= V(u)-V(\overleftarrow{u}).
\end{equation}
Then for any $a\in\r$, under $\P_a$, $(\Delta V(u);\; |u|=1)$ is distributed as $\mathcal{L}$. Let $\widetilde{Y}:=\sum_{|u|=1}\ee^{-\Delta V(u)}$ and $\widetilde{Z}:=\sum_{|u|=1}\Big(\Delta V(u)\Big)_+\ee^{-\Delta V(u)}$ so that $\P_a\Big[\Big(\widetilde{Y},\,\widetilde{Z}\Big)\in\cdot\Big]=\P[(Y,\,Z)\in\cdot]$.

Note that under $\P_{\zeta_n}$,
\begin{eqnarray}\label{eq:x}
X&=&\frac{\sum_{|u|=1}R(V(u))\ee^{-V(u)}1_{(V(u)>0)}}{R(\zeta_n)e^{-\zeta_n}}\nonumber\\
 &=&\frac{\sum_{|u|=1}R(\zeta_n+\Delta V(u))\ee^{-\Delta V(u)}1_{(\Delta V(u)>-\zeta_n)}}{R(\zeta_n)},
 \end{eqnarray}
where $(\Delta V(u);\; |u|=1)$ is independent of $\zeta_n$. By (\ref{eq:upforr}), it follows that
\begin{eqnarray*}
X&\leq & \frac{\sum_{|u|=1}c_2(\zeta_n+1)\ee^{-\Delta V(u)}1_{(\Delta V(u)>-\zeta_n)}}{R(\zeta_n)}+\frac{\sum_{|u|=1}c_2\Delta V(u)\ee^{-\Delta V(u)}1_{(\Delta V(u)>-\zeta_n)}}{R(\zeta_n)}\\
&\leq&\frac{c_2}{c_1}\sum_{|u|=1}\ee^{-\Delta V(u)}+ c_2 \frac{\sum_{|u|=1}\Delta V(u)_+\ee^{-\Delta V(u)}}{R(\zeta_n)}\\
&\leq& c_{14}\Big(\widetilde{Y}+\frac{\widetilde{Z}}{R(\zeta_n)}\Big)\leq 2c_{14}\max\Big\{\widetilde{Y},\,\frac{\widetilde{Z}}{R(\zeta_n)}\Big\},
\end{eqnarray*}
where $\Big(\widetilde{Y},\,\widetilde{Z}\Big)$ is independent of $\zeta_n$. This implies that
\begin{eqnarray}
&&\sum_{n\geq1}\E_{\zeta_n}\Big[X\Big(R(\zeta_n)\ee^{-\zeta_n}X\wedge 1\Big)\Big]\nonumber\\
&\leq &c_{15}\bigg(\sum_{n\geq0}\E\Big[\widetilde{Y}\Big(R(\zeta_n)\ee^{-\zeta_n}\widetilde{Y}\wedge 1\Big)\Big\vert \zeta_n\Big]+\sum_{n\geq0}\frac{1}{R(\zeta_n)}\E\Big[\widetilde{Z}\Big(\ee^{-\zeta_n}\widetilde{Z}\wedge 1\Big)\Big\vert\zeta_n\Big]\bigg)\nonumber\\
&=:&c_{15}\Big(\Sigma_1+\Sigma_2\Big).
\end{eqnarray}
Hence we only need to prove that
\begin{equation}\label{eq:iftrue}
\E\Big[Y\Big(\log_+Y\Big)^2\Big]+\E\Big[Z\log_+Z\Big]<\infty\Longrightarrow \E\Big[\Sigma_1\Big]+\E\Big[\Sigma_2\Big]<\infty,
\end{equation}
which leads to (\ref{eq:finimean}). On the one hand, as (\ref{eq:upforr}) gives that  $R(x)\leq c_{16}e^{x/2}$ for all $x\geq0$, we see that
\begin{eqnarray*}
\E\Big[\Sigma_1\Big]&\leq&c_{17}\E\bigg[ \sum_{n\geq0}\E\Big[\widetilde{Y}\Big(\ee^{-\zeta_n/2}\widetilde{Y}\wedge 1\Big)\Big\vert \zeta_n\Big]\bigg]\\
&=&c_{17}\sum_{n\geq0}\E\bigg[\Big(\widetilde{Y}\Big)^2\ee^{-\zeta_n}1_{\{\widetilde{Y}\leq \ee^{\zeta_n/2}\}}+\widetilde{Y}1_{\{\widetilde{Y}> \ee^{\zeta_n/2}\}}\bigg]\\
&=&c_{17}\E\bigg\{\Big(\widetilde{Y}\Big)^2\E\Big[\sum_{n\geq0}\ee^{-\zeta_n}1_{\{\zeta_n\geq 2\log \widetilde{Y}\}}\Big\vert\widetilde{Y}\Big]+\widetilde{Y}\E\Big[\sum_{n\geq0}1_{\{\zeta_n<2\log \widetilde{Y}\}}\Big\vert\widetilde{Y}\Big]\bigg\},
\end{eqnarray*}
where $\widetilde{Y}$ and $(\zeta_n)$ are independent. By (\ref{eq:expsum}),
\begin{eqnarray}
\E\Big[\Sigma_1\Big]&\leq& c_{17}\E\bigg[\Big(\widetilde{Y}\Big)^2\int_{2\log_+ \widetilde{Y}}^{\infty}\ee^{-x}R(x)U(\d x)+\widetilde{Y}\int_0^{2\log_+ \widetilde{Y}}R(x)U(\d x)\bigg],
\end{eqnarray}
which by (\ref{eq:upforr}) and (\ref{eq:upforu}) implies that
\begin{eqnarray}
\E\Big[\Sigma_1\Big]&\leq&  c_{17}c_2\E\bigg[\Big(\widetilde{Y}\Big)^2\int_{2\log_+ \widetilde{Y}}^{\infty}\ee^{-x}(x+1)U(\d x)+\widetilde{Y}\int_0^{2\log_+ \widetilde{Y}}(x+1)U(\d x)\bigg]\\
&\leq & c_{18}\E\Big[\widetilde{Y}\Big(1+\log_+\widetilde{Y}\Big)^2\Big]=c_{18}\E\Big[Y\Big(1+\log_+Y\Big)^2\Big]
\end{eqnarray}
On the other hand, in the same way, we obtain that
\begin{equation}
\E\Big[\Sigma_2\Big]\leq c_{19}\E\bigg[Z\Big(1+\log_+ Z\Big)\bigg].
\end{equation}
Consequently,
\begin{equation}
\E\Big[\Sigma_1\Big]+\E\Big[\Sigma_2\Big]\leq c_{20}\bigg(\E\Big[Y+Z\Big]+\E\Big[Y\Big(\log_+Y\Big)^2\Big]+\E\Big[Z\log_+Z\Big]\bigg).
\end{equation}
Note that (\ref{eq:settings}) ensures that $\E\Big[Y+Z\Big]<\infty$. The (\ref{eq:iftrue}) is thus proved and we completes the proof of (\ref{eq:if}).
\end{proof}

\subsection{(\ref{cond:iff}) is a necessary condition}
\label{ss:proofofonlyif}
This subsection is devoted to proving that
\begin{equation}\label{eq:onlyif}
\max\Big\{\E\Big[Z\log_+Z\Big],\; \E\Big[Y\Big(\log_+Y\Big)^2\Big]\Big\}=\infty
\Longrightarrow \E[D_\infty^{(0)}]=0.
\end{equation}

\begin{proof}[Proof of (\ref{eq:onlyif})]
According to (ii) of Theorem \ref{thm:BK}, we only need to show that
\begin{equation}\label{intermedia}
\forall y>0,\, \P\textrm{-a.s. }\sum_{n=1}^\infty\E_{\zeta_n}\Big[X;\; R(\zeta_n)\ee^{-\zeta_n}X\geq y\Big]=\infty.
\end{equation}

We break the assumption on the left-hand side of (\ref{eq:onlyif}) up into three cases. In each case, we find out a different lower bound for $X$ to establish (\ref{intermedia}). It hence follows that $D_\infty^{(0)}$ is trivial as $\E[D_\infty^{(0)}]=0$. The three cases are stated as follows:
\begin{subequations}
\begin{align}
        \E[Y(\log_+Y)^2]&=\infty,\qquad \E[Y(\log_+Y)]<\infty;\label{cond-a}\\
        \E[Y(\log_+Y)]&=\infty;\label{cond-b}\\
        \E[Z(\log_+Z)]&=\infty.\label{cond-c}
\end{align}
\end{subequations}

{\it Proof of (\ref{intermedia}) under (\ref{cond-a})}
Recall that for any particle $x\in\mathbb{T}\setminus\{\varnothing\}$, $\Delta V(u)= V(u)-V(\overleftarrow{u})$, and that under $\P_a$, $(\Delta V(u);\; |u|=1)$ is distributed as $\mathcal{L}$. For any $s\in\r$, we define a pair of random variables:
\begin{equation}
Y_+(s):=\sum_{|u|=1}\ee^{-\Delta V(u)}1_{(\Delta V(u)>-s)},\qquad Y_-(s):=\sum_{|u|=1}\ee^{-\Delta V(u)}1_{(\Delta V(u)\leq -s)}.
\end{equation}
Clearly, $\widetilde{Y}=Y_+(s)+Y_-(s)$.

It follows from (\ref{eq:x}) and (\ref{eq:upforr}) that under $\P_{\zeta_n}$,
 \begin{eqnarray*}
 X&\geq & \frac{\sum_{|u|=1}c_1(1+\zeta_n+\Delta V(u))\ee^{-\Delta V(u)}1_{(\Delta V(u)> -\zeta_n/2)}}{c_2(1+\zeta_n)},\\
&\geq& \frac{\sum_{|u|=1}c_1(1/2+\zeta_n/2)\ee^{-\Delta V(u)}1_{(\Delta V(u)> -\zeta_n/2)}}{c_2(1+\zeta_n)}\geq c_{21} Y_+(\zeta_n/2),
\end{eqnarray*}
where $\Big\{\Big(Y_+(s),\,Y_-(s)\Big);\; s\in\r\Big\}$ is independent of $\zeta_n$ and $c_{21}:=\frac{c_1}{2c_2}>0$. We thus see that the assertion that for any $y>0$,
\begin{equation}
\sum_{n=1}^\infty\E\Big[Y_+(\zeta_n/2);\; R(\zeta_n/2)\ee^{-\zeta_n}Y_+(\zeta_n/2)\geq y\Big\vert \zeta_n\Big]=\infty,\qquad \P\textrm{-a.s.,}
\end{equation}
yields (\ref{intermedia}). It is known that $\zeta_n\rightarrow\infty$ as $n$ goes to infinity (see, for example, \cite{biggins03}). It suffices that
\begin{equation}\label{eq:case1}
\sum_{n=1}^\infty F(\zeta_n/2, \zeta_n)=\infty, \qquad \P\textrm{-a.s.}
\end{equation}
where
\begin{equation}
F(s,z):=\E\Big[Y_+(s);\; \log Y_+(s)\geq z\Big],\quad s, z\in\r.
\end{equation}

Let $F_1(z):=\E[Y;\; \log Y\geq z]$ which is positive and non-increasing. It follows from Lemma \ref{many-to-one} and (\ref{eq:settings}) that $\E[Y]=1$. Therefore, for any $s$, $z\in\r$,
\begin{equation}
0\leq F(s,z)\leq F_1(z)\leq \E[Y]=1.
\end{equation}

On the one hand,  we deduce from (\ref{cond-a}) that
\begin{eqnarray*}
\int_0^\infty F_1(y)y\d y&=&\int_0^\infty \E\Big[Y1_{(\log Y\geq y)}\Big]y\d y=\E\bigg[Y\int_0^{(\log_+ Y)} y\d y;\; Y\geq 1\bigg]\\
                           &= & \E\Big[Y(\log_+ Y)^2\Big]/2=\infty.
\end{eqnarray*}
According to Proposition \ref{prop:spine}, $\P$-almost surely,
\begin{equation}\label{eq:sumone}
\sum_{n=1}^\infty F_1(\zeta_n)=\infty.
\end{equation}

On the other hand, we can prove that $\sum_{n=1}^\infty \big[F_1(\zeta_n)-F(\zeta_n/2,\zeta_n)\big]<\infty$, $\P$-a.s. In fact, as  $Y=Y_+(s)+Y_-(s)$ under $\P$, for any $s$, $y\in\r$,
\begin{eqnarray*}
F_1(y)-F(s,y)&=&\E\Big[Y1_{(\log Y\geq y)}-Y_+(s)1_{(\log Y_+(s)\geq y)}\Big]\\
           &=&\E\Big[Y1_{(\log Y\geq y>\log Y_+(s))}+Y1_{(\log Y_+(s)\geq y)}-Y_+(s)1_{(\log Y_+(s)\geq y)}\Big]\\
           &=&\E\Big[Y1_{(\log Y\geq y>\log Y_+(s))}+Y_-(s)1_{(\log Y_+(s)\geq y)}\Big].
\end{eqnarray*}
Note that $Y\leq 2\max\{Y_+(s), Y_-(s)\}$ under $\P$. It follows that
\begin{eqnarray*}
 F_1(y)-F(s,y)  &\leq &\E\Big[2Y_-(s)1_{(\log Y\geq y>\log Y_+(s),\; Y_+(s)\leq Y_-(s))}+Y1_{(\log Y\geq y>\log Y_+(s),\; Y_+(s)>Y_-(s))}\Big]\\
                        & &+\E\Big[Y_-(s)1_{(\log Y_-(s)\geq y)}\Big]\\
           &\leq &3\E\Big[Y_-(s)\Big]+\E\Big[Y1_{(\log Y\geq y>\log Y_+(s),\; Y_+(s)> Y_-(s))}\Big]\\
           &\leq &3\E\Big[Y_-(s)\Big]+\E\Big[Y1_{(\log Y\geq y>\log (Y/2)}\Big]=:d_1(s)+d_2(y).
\end{eqnarray*}
As a consequence,
\begin{equation}\label{eq:a2up}
\sum_{n=1}^\infty \big[F_1(\zeta_n)-F(\zeta_n/2,\zeta_n)\big]\leq \sum_{n\geq 0}d_1(\zeta_n/2)+\sum_{n\geq 0}d_2(\zeta_n).
\end{equation}
Taking expectation on both sides yields that
\begin{eqnarray}\label{eq:expd1}
\E\Big[\sum_{n=1}^\infty \big(F_1(\zeta_n)-F(\zeta_n/2,\zeta_n)\big)\Big]&\leq&\E\Big[\sum_{n\geq 0}d_1(\zeta_n/2)\Big]+\E\Big[\sum_{n\geq 0}d_2(\zeta_n)\Big]\nonumber\\
&=& \int_0^\infty d_1(x/2)R(x)U(\d x)+\int_0^\infty d_2(x)R(x)U(\d x),
\end{eqnarray}
where the last equality comes from (\ref{eq:expsum}).

For the first integration, we deduce from Lemma \ref{many-to-one} that
\begin{equation}
d_1(s)=3\E\Big[Y_-(s)\Big]=3\E\Big[\sum_{|x|=1}\ee^{-V(x)}1_{(V(x)\leq -s)}\Big]=3\P(-S_1\geq s).
\end{equation}
By (\ref{eq:upforr}), (\ref{eq:upforu}) and (\ref{eq:settings}),
\begin{eqnarray*}
\int_0^\infty d_1(x/2)R(x)U(\d x) &=& 3\int_0^\infty \P(-2S_1\geq x)R(x)U(\d x)\\
                                         &=& 3\E\Big[\int_0^{-2S_1} R(x)U(\d x); -2S_1\geq 0\Big]\\
                                         &\leq&\leq c_{22}\E\Big[\Big(1+(-2S_1)_+\Big)^2\Big]<\infty.
\end{eqnarray*}

For the second integration on the right-hand side of (\ref{eq:expd1}), as $d_2(y)=\E\Big[Y1_{(\log Y\geq y>\log (Y/2)}\Big]$, we use (\ref{eq:upforr}), (\ref{eq:upforu}) and (\ref{cond-a}) to obtain that
\begin{eqnarray*}
\int_0^\infty d_2(x)R(x)U(\d x)&=&  \int_0^\infty \E\Big[Y1_{(\log Y\geq x>\log (Y/2)}\Big] R(x)U(\d x)\\
&=& \E\Big[Y \int_{(\log Y-\log 2)_+}^{\log_+Y}R(x)U(\d x)\Big]\\
&\leq& c_{23}\E\Big[Y(1+\log_+Y)\Big]<\infty.
\end{eqnarray*}

Going back to (\ref{eq:expd1}), we conclude that
\begin{equation}
\E\Big[\sum_{n=1}^\infty \big(F_1(\zeta_n)-F(\zeta_n/2,\zeta_n)\big)\Big]\leq\E\Big[\sum_{n\geq 0}d_1(\zeta_n/2)\Big]+\E\Big[\sum_{n\geq 0}d_2(\zeta_n)\Big]<\infty.
\end{equation}
Therefore, $\P$-a.s.,
\begin{equation}
\sum_{n=1}^\infty \big[F_1(\zeta_n)-F(\zeta_n/2,\zeta_n)\big]<\infty,
\end{equation}
which, combined with (\ref{eq:sumone}), implies (\ref{eq:case1}). Thus (\ref{intermedia}) is proved under (\ref{cond-a}).

{\it Proof of (\ref{intermedia}) under (\ref{cond-b})} Now we suppose that $\E[Y\log_+Y]=\infty$. By (\ref{eq:upforr}), we observe that under $\P_{\zeta_n}$,
\begin{eqnarray}\label{eq:lowerboundb}
X&=&\frac{\sum_{|u|=1}R(\Delta V(u)+\zeta_n)\ee^{-\Delta V(u)}1_{(\Delta V(u)>-\zeta_n)}}{R(\zeta_n)}\nonumber\\
 &\geq &c_1\frac{Y_+(\zeta_n)}{R(\zeta_n)},
\end{eqnarray}
where $\{Y_+(s);\; s\in\r\}$ and $\zeta_n$ are independent.

To establish (\ref{intermedia}), we only need to show that for any $y\geq1$,
\begin{equation}
\sum_{n\geq 1}\E\Big[\frac{Y_+(\zeta_n)}{R(\zeta_n)};\; Y_+(\zeta_n)\geq ye^{\zeta_n}\Big\vert\zeta_n\Big]=\sum_{n\geq 1}\frac{F(\zeta_n,\log y+\zeta_n)}{R(\zeta_n)}=\infty,\qquad \P\textrm{-a.s.}
\end{equation}

 For any $y\geq1$ fixed, let
 \begin{equation}
 F_2(x):=\frac{F_1(\log y+x)}{R(x)},\; \forall x\geq 0,
 \end{equation}
  which is non-increasing as $R(x)=U^-([0,x))$ is non-decreasing and $F_1$ is non-increasing. One sees that
  \begin{equation}\label{eq:casebcut}
\sum_{n\geq1}F_2(\zeta_n)=  \sum_{n\geq 1}\frac{F(\zeta_n,\log y+\zeta_n)}{R(\zeta_n)}+\sum_{n\geq1} \frac{F_1(\log y+\zeta_n)-F(\zeta_n,\log y+\zeta_n)}{R(\zeta_n)}.
  \end{equation}

  By (\ref{eq:upforr}), $\frac{F_1(\log y+x)}{c_2(1+x)}\leq F_2(x)\leq \frac{1}{c_1}$.  It then follows from (\ref{cond-b}) that
\begin{eqnarray*}
\int_0^\infty F_2(x)x\d x&\geq &\int_0^\infty F_1(\log y+x)\frac{x}{c_2(1+x)}\d x\\
&\geq&\int_1^\infty c_{24}\E\Big[Y1_{(\log Y\geq \log y+x)}\Big] \d x\\
&\geq& c_{24}\E[Y(\log Y-\log y-1)_+]=\infty.
\end{eqnarray*}
By Proposition \ref{prop:spine},
\begin{equation}\label{eq:casebinfini}
\sum_{n\geq0} F_2(\zeta_n)=\sum_{n\geq0} \frac{F_1(\log y+\zeta_n)}{R(\zeta_n)}=\infty,\; \P\textrm{-a.s.}
\end{equation}

In view of (\ref{eq:casebcut}) and (\ref{eq:casebinfini}), it suffices to show that $\P$-a.s.,
\begin{equation}\label{eq:caseb2}
\sum_{n\geq0} \frac{F_1(\log y+\zeta_n)-F(\zeta_n,\log y+\zeta_n)}{R(\zeta_n)}<\infty.
\end{equation}

Recall that $F_1(z)-F(s,z)\leq d_1(s)+d_2(z)$. By  (\ref{eq:expsum}),
\begin{eqnarray}\label{eq:b2up}
&&\E\Big[\sum_{n\geq 0}\frac{F_1(\log y+\zeta_n)-F(\zeta_n,\log y+\zeta_n)}{R(\zeta_n)}\Big]\\
&\leq&\E\Big[\sum_{n\geq0}\frac{d_1(\zeta_n)+d_2(\log y+\zeta_n)}{R(\zeta_n)}\Big]=\int_0^\infty \Big[d_1(x)+d_2(\log y+x)\Big]U(\d x).\nonumber
\end{eqnarray}

On the one hand, recalling that $d_1(x)=3\P(-S_1\geq x)$, we deduce from (\ref{eq:upforu}) that
\begin{eqnarray}\label{eq:bd1up}
\int_0^\infty d_1(x)U(dx)&=&\int_0^\infty 3\P(-S_1\geq x)U(dx)\\
                         &=&3\E\Big[\int_0^{(-S_1)_+}U(dx)\Big]\nonumber\\
                         &\leq& 3c_4 \E\Big[1+(-S_1)_+\Big]<\infty.\nonumber
\end{eqnarray}
On the other hand, recalling that $d_2(x)=\E[Y;\; \log Y\geq x>\log Y-\log 2]$, by (\ref{eq:upforu}) again, we obtain that
\begin{eqnarray}\label{eq:bd2up}
\int_0^\infty d_2(\log y+x)U(dx)&=&\int_0^\infty \E[Y1_{(\log Y\geq \log y+x>\log Y-\log 2)}]U(dx)\\
                               &=&\E\Big[Y\int_{(\log Y-\log y-\log 2)_+}^{(\log Y-\log y)_+}U(dx)\Big]\nonumber\\
                               &\leq &c_4(1+\log 2)\E[Y]<\infty.\nonumber
\end{eqnarray}
Combined with (\ref{eq:bd1up}) and (\ref{eq:bd2up}), (\ref{eq:b2up}) becomes that
\begin{equation}
\E\Big[\sum_{n\geq 1}\frac{F_1(\log y+\zeta_n)-F(\zeta_n,\log y+\zeta_n)}{R(\zeta_n)}\Big]<\infty.
\end{equation}
We thus get (\ref{eq:caseb2}), and completes the proof of (\ref{intermedia}) given (\ref{cond-b}).

{\it Proof of (\ref{intermedia}) under (\ref{cond-c})} In this part we assume that $\E[Z\log_+ Z]=\infty$  with $Z=\sum_{|u|=1}V(u)_+\ee^{-V(x)}\geq0$. We observe that under $\P_{\zeta_n}$,
\begin{eqnarray}\label{eq:case3}
    X &\geq&  \frac{\sum_{|u|=1}R(\Delta V(u)+\zeta_n)\ee^{-\Delta V(u)}1_{(\Delta V(u)> 0)}}{R(\zeta_n)}\nonumber\\
    &\geq& \frac{c_1}{R(\zeta_n)}\widetilde{Z},
\end{eqnarray}
where $\widetilde{Z}=\sum_{|x|=1}\Big(\Delta V(x)\Big)_+\ee^{-\Delta V(x)}$ is independent of $\zeta_n$. As a consequence, for any $y>0$,
\begin{equation}
\sum_{n\geq 1}\E_{\zeta_n}\Big[X;\; R(\zeta_n)\ee^{-\zeta_n}X\geq y\Big]\geq \sum_{n\geq 1}\frac{c_1}{R(\zeta_n)}\E\Big[\widetilde{Z};\, c_1\widetilde{Z}\geq y\ee^{\zeta_n}\Big\vert\zeta_n\Big].
\end{equation}
Recall that $\widetilde{Z}$ is distributed as $Z$ under $\P$. Therefore, it is sufficient to prove that for any $y>0$,
\begin{equation}
\sum_{n\geq 1}\frac{1}{R(\zeta_n)}\E\Big[\widetilde{Z};\, \widetilde{Z}\geq y\ee^{\zeta_n}\Big\vert\zeta_n\Big]=\sum_{n\geq 1}F_3(\zeta_n)=\infty,\qquad \P\textrm{-a.s.}
\end{equation}
where
\begin{equation}
F_3(z):=\frac{\E[Z;\; \log Z\geq z+\log y]}{R(z)},\;\forall z\geq 0.
\end{equation}
Since $R$ is non-decreasing, the function $F_3$ is non-increasing. By Lemma \ref{many-to-one} and (\ref{eq:upforr}),
\begin{equation}
0\leq F_3(z)\leq\frac{\E[Z]}{R(z)}\leq\frac{\E[(S_1)_+]}{c_1}<\infty.
\end{equation}
Moreover, by (\ref{eq:upforu}) and (\ref{cond-c}),
\begin{eqnarray}
\int_0^\infty F_3(x)x \d x &\geq &\int_1^\infty c_{25}\E\Big[Z;\; \log Z-\log y\geq x\Big]\d x\\
                               &\geq& c_{25}\E\Big[Z(\log Z-\log y-1)_+\Big]=\infty.\nonumber
\end{eqnarray}
Because of Proposition \ref{prop:spine}, we obtain that for any $y>0$,
\begin{equation}
\sum_{n\geq 1}\frac{1}{R(\zeta_n)}\E\Big[\widetilde{Z};\, \widetilde{Z}\geq y\ee^{\zeta_n}\Big\vert\zeta_n\Big]=\sum_{n\geq 1}F_3(\zeta_n)=\infty,\qquad \P\textrm{-a.s.}
\end{equation}
which completes the proof of (\ref{intermedia}) under (\ref{cond-c}).

\end{proof}

\textbf{Acknowledgements}

I am grateful to my Ph. D. supervisor Prof. Zhan SHI for his advice and help. I also wish to thank my colleagues of Laboratoire des Probabilit\'es et Mod\`eles Al\'eatoires in Universit\'e Paris 6 for the enlightening discussions.

\bigskip
\bigskip
\end{document}